\documentclass[11pt]{article}
\usepackage{enumerate}
\usepackage{amsmath,amsthm,amssymb,amsfonts}
\usepackage{mathtools}
\usepackage{graphicx}
\usepackage{color}
\usepackage{subfig}
\usepackage{url}
\usepackage{hyperref}
\usepackage{tikz}
\usetikzlibrary{positioning}
\usetikzlibrary{shapes,arrows}
\usepackage{algorithmic}
\usepackage{algorithm}

\usepackage{booktabs}
\usepackage{multirow}
\usepackage[T1]{fontenc}  
\usepackage{newpxtext, newpxmath}  

\usepackage{comment}
\usepackage{geometry}
\usepackage[round]{natbib}

\newcommand{\mini}{\text{minimize}}

\newcommand{\st}{\mbox{subject to}}
\newcommand{\T}{{\top}}
\newcommand{\vt}[1]{\boldsymbol{#1}}
\newcommand{\bb}[1]{\mathbb{#1}}

\newenvironment{prob}[1]
{
\[
\noindent
\begin{array}{lll}
\text{#1}:
}
{
\end{array}
\]
}

\theoremstyle{plain}
	  \newtheorem{thm}{Theorem}

\theoremstyle{definition}
	  \newtheorem{dfn}[thm]{Definition}

\theoremstyle{remark}

\mathtoolsset{showonlyrefs=true}
\allowdisplaybreaks
\title{
A scalable mixed-integer conic optimization approach to cardinality-constrained Poisson regression with safe screening
}

\author{
  Kota Kurihara\thanks{Joint Graduate School of Mathematics for Innovation, Kyushu University, Fukuoka, Japan (\texttt{E-mail:~kurihara.kota.881@s.kyushu-u.ac.jp})}
  \and
  Yoichi Izunaga\thanks{Faculty of Economics, Kyushu University, Fukuoka, Japan (\texttt{E-mail:~izunaga@econ.kyushu-u.ac.jp})}
}

\date{\today}

\begin{document}

\maketitle
\begin{abstract}
  This paper introduces a novel approach for cardinality-constrained Poisson regression to address feature selection challenges in high-dimensional count data. We formulate the problem as a mixed-integer conic optimization, enabling the use of modern solvers for optimal solutions. To enhance computational efficiency, we develop a safe screening based on Fenchel conjugates, thereby effectively removing irrelevant features before optimization.
  
  Experiments on synthetic datasets demonstrate that our safe screening significantly reduces the problem size, leading to substantial improvements in computational time. Our approach can solve Poisson regression problems with tens of thousands of features, exceeding the scale of previous studies. This work provides a valuable tool for interpretable feature selection in high-dimensional Poisson regression.
\end{abstract}

\section{Introduction}

Poisson regression is a powerful statistical method widely used for analyzing count data, with applications in diverse fields such as predicting the number of insurance claims, modeling website traffic, and analyzing the frequency of specific events.

Given a dataset of $n$ observations $\{ (\vt{x}_{i},y_{i}) \in \bb{R}^{m} \times \bb{N} \cup \{0\} \mid i\in \{1,2,\ldots,n\} \}$,
where $\vt{x}_{i}$ is an $m$-dimensional feature vector and $y_{i}$ is the count variable to be predicted,
we assume that the random variable $Y$ follows a Poisson distribution:
\begin{align*}
  \text{Pr}(Y=\ell \mid \mu) &= \frac{\mu^{\ell}\exp(-\mu)}{\ell!}
\end{align*}
where $\mu \in \bb{R}_{+}$ is a parameter representing both the mean and variance of the distribution.
In Poisson regression, the distribution parameter $\mu$ is explained by the linear regression model
\begin{align}
  \mu = \exp(\vt{w}^{\T} \vt{x} + b),
\end{align}
where $\vt{w} \in \bb{R}^{m}$ is the vector of coefficients, and $b \in \bb{R}$ is the intercept.
Then, the average negative log-likelihood function $L(\vt{w},b)$, often referred to as the Poisson loss, is given by
\begin{align}
  L(\vt{w},b)
  &=
  \frac{1}{n}\sum_{i=1}^{n} \left( \exp(\vt{w}^{\T} \vt{x}_{i} + b)-y_{i}\left(\vt{w}^{\T} \vt{x}_{i} + b \right) +\log y_{i}!\right).
\end{align}
The goal of Poisson regression is to estimate $\vt{w}$ and $b$ by minimizing the Poisson loss:
\begin{align*}
  \mini_{}~L(\vt{w},b)~\st~(\vt{w},b) \in \bb{R}^{m+1}.
\end{align*}

A key challenge arises when applying Poisson regression to datasets with a large number of features, a situation common in genomics, where researchers might attempt to predict disease incidence based on thousands of gene expression levels obtained through RNA-Seq~\citep{AH_GB2010, FJE_BMC2022}. In such high-dimensional settings, overfitting becomes a significant concern, leading to poor out-of-sample prediction accuracy~\citep{F_JMLR2003}. For example, the high dimensionality of gene expression data, coupled with the relatively small sample sizes typical of many genomic studies, makes these models highly susceptible to overfitting, hindering their ability to generalize to new data \citep{THAI_GB2019}.

To address this issue, regularization methods such as Lasso ($\ell_{1}$-regularization) and Ridge ($\ell_{2}$-regularization) are often employed to reduce model complexity~\citep{HTT_STS2020}. However, these methods do not explicitly control the number of non-zero coefficients in the model. In contrast, cardinality-constrained regression offers a direct approach to controlling model sparsity by limiting the number of features included in the model. This is particularly valuable when interpretability is paramount, as it allows for the identification of a small subset of key features that drive the prediction.

In this paper, we consider the Poisson regression with $\ell_{2}$-regularization and cardinality constraint, which is formulated as follows:
\begin{align}
  \displaystyle \mini_{(\vt{w},b)\in \bb{R}^{m+1}} L(\vt{w},b) + \frac{1}{\gamma} \|\vt{w}\|_{2}^{2} ~\st~ \|\vt{w}\|_{0} \leq k,
\end{align}
where the $\ell_{0}$-pseudo norm $\|\vt{w}\|_{0}$ of a vector $\vt{w}$ counts the number of nonzeros in $\vt{w}$, 
$\gamma > 0$ is a regularization parameter and $k$ is a user-defined parameter for specifying the number of non-zero coefficients.
This optimization problem presents two challenges: (i) the nonlinearity induced by the exponential term in the Poisson loss function and (ii) the combinatorial and nonlinear natures arising from the cardinality constraint. Indeed, the cardinality constraint makes the problem NP-hard~\citep{N_SICOMP1995}.

\paragraph{Related work}

Feature selection, the process of identifying the most relevant variables for a predictive model, is a fundamental challenge in statistical modeling and machine learning.
The most straightforward approach, evaluating the loss function for every possible combination of features, quickly becomes infeasible.
The number of combinations grows exponentially, rendering an exhaustive search impractical for even moderately sized datasets.
Faced with this computational hurdle, researchers have developed a variety of alternative strategies. Stepwise methods~\citep{E_MMDC1960, A_2002},
for example, offer a more efficient approach by iteratively adding or removing features based on their individual contribution to the model's performance.
While popular, stepwise methods are inherently greedy and do not guarantee finding the globally optimal feature subset.
Metaheuristic algorithms, such as those inspired by nature~\citep{A_SOIC_2019} or particle swarm behavior~\citep{AK_CSTM2018}, have also been explored for Poisson regression feature selection.
However, these methods are also heuristic in nature and lack the ability to verify the optimality of their solutions.

In recent years, mixed-integer optimization (MIO) has emerged as a powerful tool for feature selection. Driven by advances in both algorithms and computing hardware, MIO offers the distinct advantage of providing optimal solutions with respect to various feature selection criteria~\citep{MT_EJOR2015, GP_INFORMS2021}.
For a comprehensive overview of MIO approaches to feature selection, see \cite{BKM_AOR2016} and \cite{DHHM_JMLR2021}.
\cite{SKT_PLOSONE2021} further demonstrated the potential of MIO by proposing a mixed-integer quadratic optimization (MIQO) formulation for sparse Poisson regression.
By employing piecewise-linear approximations to the Poisson loss function, they achieved improved out-of-sample prediction performance compared to existing techniques.

Despite these advances, existing MIO approaches often struggle with the computational demands of large, high-dimensional datasets.
The experiments conducted by~\cite{SKT_PLOSONE2021},
for instance, were limited to problems with fewer than 40 features.
This limitation highlights the need for more efficient methods capable of tackling the challenges of cardinality-constrained Poisson regression in modern, large-scale applications.

To address this need, researchers have explored techniques to strengthen the relaxations used within MIO solvers.
A common strategy is to apply a branch-and-bound method, where tighter bounds from the continuous relaxation of the underlying problem lead to a more efficient search.
\cite{GL_MP2010} proposed perspective reformulations to derive stronger continuous relaxations.

Several studies of MIO approaches to feature selection applied linear approximations to the nonlinear loss functions such as the Poisson loss~\citep{SKT_PLOSONE2021} and the logistic loss~\citep{STMY_COAP2016}. On the other hand, there has been a growing interest in non-symmetric cone optimization, which extends the traditional symmetric cone optimization framework including linear programming, second-order cone programming, and semidefinite programming. Non-symmetric cones, such as the exponential cone and the power cone, have been shown to be powerful tools for modeling a wide range of optimization problems (see, e.g.,\citep{K_IJOC2021} and \citep{YX_arxiv2021}).
In the literature on feature selection, \cite{AK_2023} proposed a mixed-integer exponential cone formulation for sparse logistic regression.

In this paper,
we propose a mixed-integer conic optimization formulation for cardinality-constrained Poisson regression, enhanced by a safe screening originally developed by~\citet{AG_ICML2020} for ordinary least squares regression. 
The safe screening efficiently identifies key features and removes irrelevant features before the optimization process, significantly reducing the problem size and improving computational efficiency.

\paragraph{Contributions}
The main contributions of this paper are as follows:

\begin{enumerate}
  \item We propose a novel mixed-integer conic optimization formulation for cardinality-constrained Poisson regression, leveraging the exponential cone representation and the perspective function of the $\ell_{2}$-regularization term to derive a stronger continuous relaxation.
  \item We derive a safe screening rule for the Poisson loss minimization problem by exploiting the Fenchel conjugate function of the perspective function of the $\ell_{2}$-regularization term, enabling the efficient identification of key features and removal of irrelevant features.
  \item We develop a greedy heuristic, based on the optimal dual solution, to efficiently compute an upper bound required for the safe screening. As a valuable byproduct, this heuristic also provides a high-quality initial solution for the Branch-and-Bound process.
  \item We conduct computational experiments on synthetic datasets, demonstrating that our safe screening, enabled by the greedy heuristic for upper bound estimation, significantly reduces problem size and improves the computational efficiency of the mixed-integer conic optimization approach. Notably, our approach can solve Poisson regression problems with a number of features in the order of ten-thousand, which is significantly larger than the size of datasets used in previous studies. The results show that the safe screening effectively identifies the best subset of features in many instances and reduces the optimality gap compared to other methods. The greedy heuristic contributes to quickly obtaining high-quality solutions, even when the safe screening is not effective, by providing a good starting point for the optimization process.
\end{enumerate}

The remainder of this paper is organized as follows: In Section 2, we provide preliminaries on convex cones, perspective functions, and Fenchel conjugate functions, and we reformulate the Poisson loss minimization problem as a mixed-integer conic optimization problem. In Section 3, we derive the safe screening rule. In Section 4, we present our computational experiments. Finally, we conclude the paper in Section 5.

\section{Perspective reformulations and conic representations}

\subsection{Notation and preliminaries}

Throughout this paper, we let nonbold faced characters such as $x$ denote scalars, bold faced characters such as $\vt{x}$ denote vectors,
and $\vt{0}$ denote a vector of all 0's, respectively.
We denote the set of consecutive positive integers as $[n]=\{ 1, 2,\ldots, n\}$, the set of $n$-dimensional real vectors as $\bb{R}^{n}$, and the set of $n$-dimensional non-negative real vectors as $\bb{R}_{+}^{n}$.

We first introduce some notions, namely two non-linear cones, perspective function, and Fenchel conjugate function, which are used in later discussions.

\begin{dfn}
  The exponential cone is defined as
  \begin{align}
    K_{\exp} = \left\{ (x_{1},x_{2},x_{3}) \in \bb{R}^{3} \mid x_{1} \ge x_{2}\exp{(x_{3}/x_{2}), x_{2} >0} \right\} \cup \{ (x_{1},0,x_{3}) \in \bb{R}^{3} \mid x_{1}\ge 0, x_{3}\le 0 \}. \label{eq:exp_cone}
  \end{align}
\end{dfn}

\begin{dfn}
  The three-dimensional rotated quadratic cone is defined as
  \begin{align}
    K_{\text{rq}}^{3} = \left\{ (x_{1},x_{2},x_{3}) \in \bb{R}^{3} \mid 2x_{1}x_{2} \ge x_{3}^{2}, x_{1} \ge 0, x_{2} \ge 0 \right\}. \label{eq:rq_cone}
  \end{align}
\end{dfn}

The above two cones are known to be convex, but the exponential cone is not symmetric unlike the rotated quadratic cone.
The exponential cone can be used to model a variety of constraints involving exponentials and logarithms. See~\cite{Mosek_cookbook} for more details.

\begin{dfn}
  The perspective function of a given function $f:\bb{R}^{d} \to \bb{R}$ is the function $\tilde{f}: \bb{R}^{d+1} \to \bb{R}\cup \{+\infty\}$ defined as follows:
  \begin{align}
    \tilde{f}(\lambda,\vt{x}) &=
    \begin{cases}
      \lambda f(\vt{x}/\lambda) & \text{if $\lambda>0$},\\
      0 & \text{if $\lambda=0$ and $\vt{x}=\vt{0}$},\\
      +\infty & \text{otherwise}.
    \end{cases}
  \end{align}
\end{dfn}

For a convex function $f$, the perspective function $\tilde{f}$ is also convex~\citep{HU_2011}.

\begin{dfn}
  The Fenchel conjugate function of a given function $f:\bb{R}^{d} \to \bb{R}$ is the function $f^{*}:\bb{R}^{d} \to \bb{R}$ defined as follows:
  \begin{align}
    f^{*}(\vt{y}) &= \max \left\{ \vt{y}^{\T}\vt{x} - f(\vt{x}) \mid \vt{x} \in \bb{R}^{d} \right\}. \label{eq:fenchel}
  \end{align}
\end{dfn}
From the definition~\eqref{eq:fenchel} of the Fenchel conjugate function,
we easily have
\begin{align}
  f(\vt{x}) \ge \vt{y}^{\T}\vt{x} - f^{*}(\vt{y}) \text{ for all $(\vt{x},\vt{y}) \in \bb{R}^{2d}$} \label{eq:Fenchel_Ineq}.
\end{align}

\subsection{Perspective reformulations}

Let $\vt{z}=(z_{j})_{j \in [m]} \in \{0,1\}^{m}$ be a vector of binary variables to indicate a feature is selected ($z_{j}=1$) or not ($z_{j}=0$).
We can state Poisson regression with $\ell_{2}$-regularization and cardinality constraint as the following optimization problem
\begin{prob}{(Q)}
  &\mini_{} & \displaystyle L(\vt{w},b)+\frac{1}{\gamma} \sum_{j=1}^{m}w_{j}^{2}\\
  &\st & \displaystyle \sum_{j=1}^{m} z_{j} \le k,\\
  & & \displaystyle w_{j} = w_{j}z_{j} \quad \text{for $j \in [m]$},\\
  & & \displaystyle (\vt{w},b) \in \bb{R}^{m+1}, \vt{z} \in \{0,1\}^{m}.
\end{prob}

The nonlinear constraint $w_{j} = w_{j}z_{j}$ can be linearized by
\begin{align}
  -Mz_{j} \le w_{j} \le Mz_{j}, \label{eq:bigM}
\end{align}
where $M \in \bb{R}_{+}$ is a sufficiently large constant, which is called big-M technique.
Note that $M$ is not known a priori, and a small value of $M$ may lead to a solution different from the original problem.
On the other hand, a large value of $M$ may lead to
weak continuous relaxations, and the numerical stability of the solution may be compromised.

To derive a stronger continuous relaxation, we exploit the perspective function of the $\ell_{2}$-regularization.
Replacing the $\ell_{2}$-regularization term $w_{j}^{2}$ in the objective function with its perspective function $w_{j}^{2}/z_{j}$ and remove the nonlinear constraint $w_{j} = w_{j}z_{j}$, we obtain an exact reformulation of the original problem (Q) as follows
\begin{prob}{(P)}
  &\mini_{} & \displaystyle L(\vt{w},b)+\frac{1}{\gamma} \sum_{j=1}^{m}\frac{w_{j}^{2}}{z_{j}}\\
  &\st & \displaystyle \sum_{j=1}^{m} z_{j} \le k,\\
  & & \displaystyle (\vt{w},b) \in \bb{R}^{m+1}, \vt{z} \in \{0,1\}^{m}.
\end{prob}

\subsection{Conic representation of Poisson regression}

We reformulate the Poisson loss minimization problem ($\min_{(\vt{w},b)} L(\vt{w},b))$ as a conic optimization problem by using the exponential cone representation.
Introducing an auxiliary variable $t_{i} \in \bb{R}$, the epigraph of the exponential term $t_{i} \ge \exp(\vt{w}^{\top}\vt{x}_{i}+b)$ in the loss function is represented as
\begin{align}
  (t_{i}, 1, \vt{w}^{\top}\vt{x}_{i}+b) \in K_{\exp}. \label{eq:conic_rep1}
\end{align}

Moreover, introducing non-negative auxiliary variables $s_{j} \in \bb{R}_{+}$ where  $s_{j}$ models $w_{j}^{2}$ and add the constraint $w_{j}^{2} \le s_{j}z_{j}$, which is rotated quadratic cone representable.
Namely, it holds that
\begin{align}
  w_{j}^{2} \le s_{j}z_{j} \quad \Leftrightarrow \quad (s_{j}/2, z_{j}, w_{j}) \in K_{\text{rq}}^{3}. \label{eq:conic_rep2}
\end{align}
These lead to the following reformulation of the problem (P),
which is a mixed-integer conic optimization problem
\begin{prob}{($\text{P}_{\text{conic}}$)}
  &\mini_{} & \displaystyle \frac{1}{n}\sum_{i=1}^{n} (t_{i}-y_{i}(\vt{w}^{\top}\vt{x}_{i}+b)+\log{y_{i}!}) +\frac{1}{\gamma} \sum_{j=1}^{m}s_{j}\\
  &\st & \displaystyle \sum_{j=1}^{m} z_{j} \le k,\\
  & & \displaystyle(t_{i}, 1, \vt{w}^{\top}\vt{x}_{i}+b) \in K_{\exp} \quad \text{for $i \in [n]$},\\
  & & \displaystyle (s_{j}/2, z_{j}, w_{j}) \in K_{\text{rq}}^{3} \quad \text{for $j \in [m]$},\\
  & & \displaystyle (\vt{w},b) \in \bb{R}^{m+1}, \vt{z} \in \{0,1\}^{m}, (\vt{s},\vt{t}) \in \bb{R}_{+}^{m+n}.
\end{prob}

\section{Safe screening}

The problem (P), or equivalently the problem ($\text{P}_{\text{conic}}$), is a mixed-integer optimization problem, which is generally hard to solve.
If we can identify a subset of features that are guaranteed to be selected or not at any optimal solution, we can reduce the problem size,
and  an exponential speedup can be expected with respect to the number of fixed variables in enumerative methods such as the branch-and-bound.
This problem reduction method is called \textit{safe screening} in machine learning community~\citep{BERG_EUSIPCO2014}. In operations research community, it is also called \textit{pegging test},
and favorable results have been reported on Knapsack problem~\cite{N_MS1976}, set-covering problem~\citep{BC_OR1996} and so on.

We derive the safe screening rule for the problem (P) by exploiting the Fenchel conjugate function of the perspective function of the $\ell_{2}$-regularization term proposed by~\cite{AG_ICML2020}.

\paragraph{Relaxation problems}

We denote (RP) as the relaxation problem of (P) by replacing the binary variable $z_{j} \in \{0,1\}$ with a continuous variable $z_{j} \in [0,1]$.

In particular, we consider the Fenchel conjugate function $p^{\star}$ of the perspective function $p(w_{j},z_{j})=w_{j}^{2}/z_{j}$:
\begin{align}
  p^{\star}(\lambda_{j},\zeta_{j})=\max \left\{ \lambda w_{j}+\zeta_{j} z_{j} - \frac{w_{j}^{2}}{z_{j}} \mid (w_{j},z_{j}) \in \bb{R}^{2} \right\} \quad \text{for any $\lambda_{j}, \zeta_{j} \in \bb{R}$}.
\end{align}
From Fenchel Inequality~\eqref{eq:Fenchel_Ineq}, we have
\begin{align}
  \frac{w_{j}^{2}}{z_{j}} \ge \lambda w_{j} + \zeta_{j} z_{j} - p^{\star}(\lambda_{j},\zeta_{j}) \quad \text{for any $(w_{j},z_{j},\lambda_{j},\zeta_{j}) \in \bb{R}^{4}$}. \label{eq:Fenchel_Ineq2}
\end{align}
Therefore, replacing the perspective terms with the right-hand side of the above Inequality~\eqref{eq:Fenchel_Ineq2}, we obtain the following lower envelope of the objective function of (P):
\begin{align}
  L(\vt{w},b)+\frac{1}{\gamma}\sum_{j=1}^{m}\frac{w_{j}^{2}}{z_{j}} \ge L(\vt{w},b)+\frac{1}{\gamma}\sum_{j=1}^{m}\left( \lambda_{j} w_{j} + \zeta_{j} z_{j} - p^{\star}(\lambda_{j},\zeta_{j}) \right)
\end{align}
for any multiplier vectors $(\vt{\lambda},\vt{\zeta}) \in \bb{R}^{2m}$.
Moreover, the closed form of $p^{\star}$ is provided in \cite{AG_ICML2020} as follows:
\begin{align}
  p^{\star}(\lambda_{j},\zeta_{j}) =
  \begin{cases}
    0 & \text{if $\zeta_{j}=-\lambda_{j}^{2}/4$},\\
    \infty & \text{otherwise}.
  \end{cases}
\end{align}
Then, for any single multiplier vector $\vt{\lambda} \in \bb{R}^{m}$, we have the following relaxation problem
\begin{prob}{(RP($\vt{\lambda}$))}
  &\mini_{} & \displaystyle L(\vt{w},b)+\frac{1}{\gamma} \sum_{j=1}^{m} \lambda_{j}w_{j}-\sum_{j=1}^{m}\left( \frac{\lambda_{j}^{2}}{4\gamma} \right) z_{j}\\
  &\st & \displaystyle \sum_{j=1}^{m} z_{j} \le k,\\
  & & \displaystyle (\vt{w},b) \in \bb{R}^{m+1}, \vt{z} \in [0,1]^{m}.
\end{prob}

Given a multiplier vector $\vt{\lambda} \in \bb{R}^{m}$,
we can easily solve the problem (RP($\vt{\lambda}$)) due to the decomposed structure of the problem, which separates the variables $(\vt{w},b)$ and $\vt{z}$.
Specifically, the optimization with respect to $(\vt{w},b)$ is reduced to finding a stationary point of the smooth function $L(\vt{w},b)+\frac{1}{\gamma} \sum_{j=1}^{m} \lambda_{j}w_{j}$,
which can be done by solving the following system of equations:
\begin{align}
  \nabla_{\vt{w}}L(\vt{w},b)+\frac{1}{\gamma}\vt{\lambda}=\vt{0} \quad \Leftrightarrow \quad
  \frac{\gamma}{n}\sum_{i=1}^{n} \left( y_{i}-\exp(\vt{w}^{\top}\vt{x}_{i}+b) \right)\vt{x}_{i}=\vt{\lambda}. \label{eq:opt_w}
\end{align}
On the other hand, the optimization with respect to $\vt{z}$ is reduced to sorting the values of $(\lambda_{1}^{2},\lambda_{2}^{2},\ldots,\lambda_{m}^{2})$ in descending order and selecting the top-$k$ elements of $\vt{z}$.
Namely, we can obtain an optimal solution by simply setting as
\begin{align}
  z_{j} &=
  \begin{cases}
    1 & \text{if $\lambda_{j}^{2} \ge \theta_{k}$},\\
    0 & \text{otherwise},
  \end{cases} \label{eq:opt_z}
\end{align}
where $\theta_k$ is the $k$-th largest value of $(\lambda_{1}^{2},\lambda_{2}^{2},\ldots,\lambda_{m}^{2})$.

\paragraph{Dual problem}

The optimal value of the problem (RP($\vt{\lambda}$)) provides a lower bound on the optimal value of the original problem (P),
and it depends on a choice of the multiplier vector $\vt{\lambda}$.
The problem of finding the tightest lower bound is called the Fenchel dual problem, which is defined as follows
\begin{prob}{(FD)}
  &\displaystyle \max_{\vt{\lambda} \in \bb{R}^{m}} & \displaystyle \min_{(\vt{w},b) \in \bb{R}^{m+1}, \vt{z} \in [0,1]^{m}} \left\{ L(\vt{w},b)+\frac{1}{\gamma} \sum_{j=1}^{m} \lambda_{j}w_{j}-\sum_{j=1}^{m}\left( \frac{\lambda_{j}^{2}}{4\gamma} \right) z_{j} \mid \sum_{j=1}^{m} z_{j} \le k \right\}
\end{prob}

Henceforth, we refer to the multiplier vector $\vt{\lambda}$ as dual variable.
By strong duality theorem on convex optimization, we can interchange the order of the maximization and minimization in the problem (FD).
Namely, the optimal value of the problem (FD) is equal to that of the following interchanged problem
\begin{prob}{($\overline{\text{FD}}$)}
  &\displaystyle \min_{(\vt{w},b) \in \bb{R}^{m+1}, \vt{z} \in [0,1]^{m}} & \displaystyle \max_{\vt{\lambda} \in \bb{R}^{m}} \left\{ L(\vt{w},b)+\frac{1}{\gamma} \sum_{j=1}^{m} \lambda_{j}w_{j}-\sum_{j=1}^{m}\left( \frac{\lambda_{j}^{2}}{4\gamma} \right) z_{j} \right\}\\
  & \st & \displaystyle \sum_{j=1}^{m} z_{j} \le k.
\end{prob}

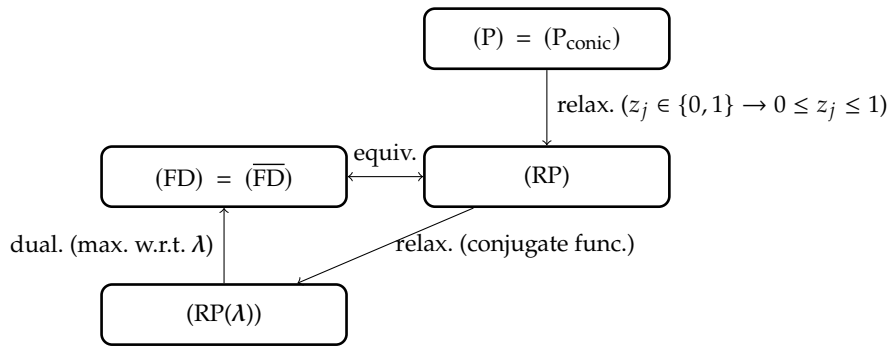
\begin{figure}[h]
\centering
\begin{tikzpicture}
  \tikzset{block/.style={rectangle, draw=black, fill=none, line width=1pt, text width=3cm, text centered, rounded corners, minimum height=0.8cm}};
  \node[block, font=\footnotesize] (P) {$(\text{P})=(\text{P}_{\text{conic}})$};
  \node[block, below=of P, font=\footnotesize] (RP) {(RP)};
  \node[block, below left=of RP, font=\footnotesize] (RRP) {(RP($\vt{\lambda}$))};
  \node[block, left=of RP, font=\footnotesize] (RD) {$(\text{FD})=(\overline{\text{FD}})$};
  
  \draw[->] (P) -- (RP) node[midway, right, font=\footnotesize] {relax.~($z_{j} \in \{0,1\} \to 0 \le z_{j} \le 1$)};
  \draw[->] (RP) -- (RRP) node[midway, right, font=\footnotesize] {relax.~(conjugate func.)};
  \draw[->] (RRP) -- (RD) node[midway, left, font=\footnotesize] {dual.~(max.~w.r.t.~$\boldsymbol{\lambda}$)};
  \draw[<->] (RP) -- (RD) node[midway, above, font=\footnotesize] {equiv.};
\end{tikzpicture}
\caption{Relationship between the problems.}
\label{fig:relationship}
\end{figure}

Note that, on the inner maximization problem of the problem ($\overline{\text{FD}}$), an optimal $\lambda_{j}$ is given by $\lambda_{j}=2w_{j}/z_{j}$ for each $j \in [m]$.
Plugging $\lambda_{j}=2w_{j}/z_{j}$ into the objective function of ($\overline{\text{FD}}$),
we obtain the same objective function of (RP).
Therefore, the problems (RP), (FD), and ($\overline{\text{FD}}$) provide the same optimal value, and are equivalent in that sense.
We summarize the relationship between the problems as mentioned above in Figure~\ref{fig:relationship}.

\paragraph{Safe screening rule}
We are ready to derive the safe screening rule for the problem (P).
Hereafter, we denote the optimal value of problem, say (P), with the additional constraint $z_{i}=\alpha$ by $v(\text{P} \mid z_{i}=\alpha)$.
Since the problem (RP($\vt{\lambda}$)) is a relaxation problem for (P), it holds that
\begin{align}
  v(\text{RP}(\vt{\lambda}) \mid z_{i}=1-\alpha)\le v(\text{P} \mid z_{i}=1-\alpha) \label{eq:peg_1}
\end{align}
for $\alpha \in \{0,1\}$.
Let $ub$ be an upper bound on the optimal value of the problem (P), 
and we suppose that
\begin{align}
  ub < v(\text{RP}(\vt{\lambda}) \mid z_{i}=1-\alpha). \label{eq:peg_2}
\end{align}
From \eqref{eq:peg_1} and \eqref{eq:peg_2}, we have $ub < v(\text{P} \mid z_{i}=1-\alpha)$.
Therefore, any feasible solution which satisfies $z_{i}=1-\alpha$ has the objective value greater than $ub$,
this implies that any optimal solutions of (P) satisfy $z_{i}=\alpha$.

Now, we focus on how to compute the objective value of the problem with a variable $z_{i}$ fixed: $v(\text{RP}(\vt{\lambda}) \mid z_{i}=1-\alpha)$ efficiently.
From \eqref{eq:opt_z}, for an index $j \in [m]$ such that $\lambda_{j}^{2} \ge \theta_{k}$, it must be $z_{j}=1$.
Owing to the problem structure and the cardinality constraint,
we can compute the optimal value of the problem (RP($\vt{\lambda}$)) with the additional constraint $z_{j}=0$ as
\begin{align}
  v(\text{RP}(\vt{\lambda}) \mid z_{j}=0)=v(\text{RP}(\vt{\lambda}))+\frac{\lambda_{j}^{2}-\theta_{k+1}}{4\gamma}. \label{eq:screening_1}
\end{align}
Similarly, for an index $j \in [m]$ such that $\lambda_{j}^{2} \le \theta_{k+1}$, it must be $z_{j}=0$, hence we readily have
\begin{align}
  v(\text{RP}(\vt{\lambda}) \mid z_{j}=1)=v(\text{RP}(\vt{\lambda}))+\frac{\theta_{k}-\lambda_{j}^{2}}{4\gamma}. \label{eq:screening_2}
\end{align}

From the above discussion, the following safe screening rules are obtained.

\begin{thm}\label{thm:ssr}
  Let $\vt{\lambda}\in \bb{R}^{m}$ be a dual variable, and $ub$ be an upper bound on the problem (P).
  Then, at any optimal solution of (P), it holds that
  \begin{align}
    z_{j}&=
    \begin{cases}
      0 & \text{if $\lambda_{j}^{2} \ge \theta_{k}$ and $v(\text{RP}(\vt{\lambda}))+(\lambda_{j}^{2}-\theta_{k+1})/4\gamma > ub$},\\
      1 & \text{if $\lambda_{j}^{2}\le \theta_{k+1}$ and $v(\text{RP}(\vt{\lambda}))+(\theta_{k}-\lambda_{j}^{2})/4\gamma > ub$.}
    \end{cases} \label{eq:ssr}
  \end{align}
\end{thm}

\paragraph{Overall description}

Safe screening rules stated by Theorem~\ref{thm:ssr} holds for any dual variable $\vt{\lambda} \in \bb{R}^{m}$.
We are interested in well-organized setting of dual variable $\vt{\lambda}$ for the effective execution of safe screening.
From the strong duality as mentioned above,
solving the problem (RP) gives us the tightest lower bound $v(\text{RP})=v(\text{RP}(\bar{\vt{\lambda}}))$,
and we can recover an optimal dual variable $\bar{\vt{\lambda}}$ by the optimality condition~\eqref{eq:opt_w}.
Indeed, let $(\vt{w}^{*},b^{*})$ be an optimal solution of (RP), the optimal $\bar{\vt{\lambda}}$ is given by
\begin{align}
  \bar{\vt{\lambda}}=\frac{\gamma}{n}\sum_{i=1}^{n} (y_{i}-\exp(\vt{x}_{i}^{\top}\vt{w}^{*}+b^{*}))\vt{x}_{i}. \label{eq:recovery}
\end{align}

In computing an upper bound $ub$ on the optimal value of the problem (P),
there exists extensive work on heuristics for sparse regression problems, such as LASSO and stepwise method.
We propose a simple greedy heuristics based on the optimal dual solution $\bar{\vt{\lambda}}$.
In particular, we select $k$ features corresponding to the top $k$-largest values of $(\bar{\lambda}_{1}^{2},\bar{\lambda}_{2}^{2},\ldots,\bar{\lambda}_{m}^{2})$,
and denote the selected features as $J_{k}$, that is $J_{k}=\{ j \in [m] \mid \bar{\lambda}_{j}^{2} \ge \theta_{k} \}$.
Then, we solve the following problem with only the selected features:
\begin{align}
  \displaystyle \min_{\vt{w}_{|J_{k}} \in \bb{R}^{k}} \, \frac{1}{n}\sum_{i=1}^{n} \left( \exp( \vt{w}_{|J_{k}}^{\top} \vt{x}_{i|J_{k}}+b)-y_{i}(\vt{w}_{|J_{k}}^{\top} \vt{x}_{i|J_{k}}+b) + \log{y_{i}!}\right)+\frac{1}{\gamma} \|\vt{w}_{|J_{k}}\|_{2}^{2} \label{eq:greedy}
\end{align}
where $\vt{x}_{i|J_{k}}=(x_{ij})_{j \in J_{k}}$ and $\vt{w}_{|J_{k}}=(w_{j})_{j \in J_{k}}$, respectively.
Any solution $(\vt{w}_{|J_{k}},b)$ of the above problem is a feasible solution of the original problem since the number of nonzero elements of $\vt{w}_{|J_{k}}$ is at most $k$, hence the objective value of the problem is a valid upper bound on the optimal value $v(P)$.
This problem is an unconstrained convex minimization problem with only $k$ features, and the objective function is smooth and strongly convex.
Therefore, it can be solved efficiently by applying first-order methods such as the gradient descent algorithm.
The overall procedure of the safe screening is summarized as follows.
\begin{algorithm}[H]
  \caption{Safe screening}\label{alg:ssr}
  \begin{algorithmic}[1]
  \STATE{Let $I_{0}$ (resp., $I_{1}$) the set of indices correspoding to the variables fixed to zero (resp., one), and initialize $I_{0}=\emptyset$, $I_{1}=\emptyset$.}
  \STATE{Solve (RP) to obtain an optimal solution $(\vt{w}^{*},b^{*})$ and the optimal value $v(\text{RP})$}
  \STATE{Compute an optimal dual variable $\bar{\vt{\lambda}}$ by Equation~\eqref{eq:recovery}}.
  \STATE{Sort the values of $(\bar{\lambda}_{1}^{2},\bar{\lambda}_{2}^{2},\ldots,\bar{\lambda}_{m}^{2})$ in descending order, and let $\theta_{\ell}$ be the $\ell$-largest value of $(\bar{\lambda}_{1}^{2},\bar{\lambda}_{2}^{2},\ldots,\bar{\lambda}_{m}^{2})$}.
  \STATE{Compute an upper bound $ub$ on the optimal value $v(\text{P})$ by solving the problem~\eqref{eq:greedy}.}
  \FOR{$j=1$ to $m$}
  \IF{$\bar{\lambda}_{j}^{2} \ge \theta_{k}$ and $v(\text{RP})+(\bar{\lambda}_{j}^{2}-\theta_{k+1})/4\gamma > ub$}
  \STATE{$I_{0}:=I_{0} \cup \{ j \}$}
  \ELSIF{$\bar{\lambda}_{j}^{2} \le \theta_{k+1}$ and $v(\text{RP})+(\theta_{k}-\bar{\lambda}_{j}^{2})/4\gamma > ub$}
  \STATE{$I_{1}:=I_{1} \cup \{ j \}$}
  \ENDIF
  \ENDFOR
  \end{algorithmic}
\end{algorithm}

Line~2 in the Algorithm~\ref{alg:ssr} is the most computationally expensive step, which requires solving the problem (RP).
However, the problem (RP), or equivalently $(\text{RP}_{\text{conic}})$, is a continuous conic optimization problem, and we can solve it by off-the-shelf solvers such as MOSEK~\citep{Mosek_11}.

\section{Numerical experiments}

This section presents numerical results evaluating the effectiveness of the safe screening proposed in the previous section, using synthetic datasets. Our experiments focus solely on the impact of the safe screening on problem size reduction and computational time, as out-of-sample prediction performance has already been investigated by~\cite{SKT_PLOSONE2021}.

All experiments were conducted on a Mac Studio equipped with a 24-core Apple M2 Ultra processor (3.50~GHz) and 192GB of memory. All methods were implemented in Julia 1.11, and MOSEK 11 was used to solve the conic optimization problems.

\subsection{Data generation}

Our experiments consider the high-dimensional case ($m > n$), which is common in genomics. Following the experimental setup in prior studies~\citep{BKM_AOR2016, HTT_STS2020, SKT_PLOSONE2021}, we generate synthetic datasets using the following parameters: the number of candidate features was set to $m=10000$, the true sparsity to $k_{\text{true}}=30$, the maximum value of the count variable to $y_{\text{max}}=10$, and the number of observations to $n \in \{1000, 2000\}$.

First, a set $S^{*}$ of $k_{\text{true}}$ features is drawn uniformly from the $m$ candidate features,
and define a vector of the true coefficients $\vt{w}^\text{true}$ as follows:
\begin{align}
  w_{j}^\text{true} =
  \begin{cases}
    1 & \text{if $j \in S^{*}$},\\
    0 & \text{otherwise}.
  \end{cases}
\end{align}

Then, we sample vectors $\vt{x}_{i} \in \bb{R}^{m}$ from an $m$-dimensional multivariate standard normal distribution as $\vt{x}_{i} \sim N(\vt{0, \Sigma})$, where $\Sigma \in \bb{R}^{m \times m}$ is a covariance matrix whose $(i,j)$-th element is $\rho^{|i-j|}$. The parameter $\rho$ represents the correlation strength between features.
We also sample the error term form a normal distribution as $\varepsilon_{i} \sim N(0,\sigma^{2})$, where $\sigma^{2}$ is the variance.
Finally, we generate the count variable $y_{i}$ by rounding
\begin{align*}
  \exp\left( \frac{(\vt{w}^{\text{true}})^{\top}\vt{x}_{i}}{\sqrt{(\vt{w}^{\text{true}})^{\top}\Sigma \vt{w}^{\text{true}}}} +\varepsilon_{i} \right)
\end{align*}
to the nearest integer.
When the generated count variable $y_{i}$ is larger than $y_\text{max}$, we set it to $y_\text{max}$.
We compare the performance of the safe screening in six different regimes of noise level ($\sigma^{2}$) and correlation strength ($\rho$) described in Table~\ref{tab:settings}, and conducted five trials for each setting to evaluate the consistency and robustness of the results.

\begin{table}[h]
  \centering
  \caption{Experimental regimes.}
  \label{tab:settings}
  \begin{tabular}{l c c}
    \toprule
    & Low correlation & High correlation \\
    \midrule
    Low noise & $(\sigma^{2}, \rho)=(0.01, 0.35)$ & $(\sigma^{2}, \rho)=(0.01, 0.70)$ \\
    Medium noise & $(\sigma^{2}, \rho)=(0.10, 0.35)$ & $(\sigma^{2}, \rho)=(0.10, 0.70)$ \\
    High noise & $(\sigma^{2}, \rho)=(1.00, 0.35)$ & $(\sigma^{2}, \rho)=(1.00, 0.70)$ \\
    \bottomrule
  \end{tabular}
\end{table}

\subsection{Performance of safe screening}

First, we examine the effect of the safe screening on problem size reduction.
\cite{BPP_2020} suggest that setting the regularization parameter $\gamma=1/\sqrt{n}$ is an appropriate scaling for regression problems.
To better understand the effect of the safe screening, we set the regularization parameter $\gamma$ to $\gamma_{0}=1/\sqrt{n}$ and varied it by multiplying $\gamma_{0}$ by $2^{\ell}$ for $\ell \in \{0, 2, 4\}$. The parameter $k$ in the cardinality constraint is set to $k=30$.

\begin{figure}[h]
  \centering
  \includegraphics[width=1.0\textwidth]{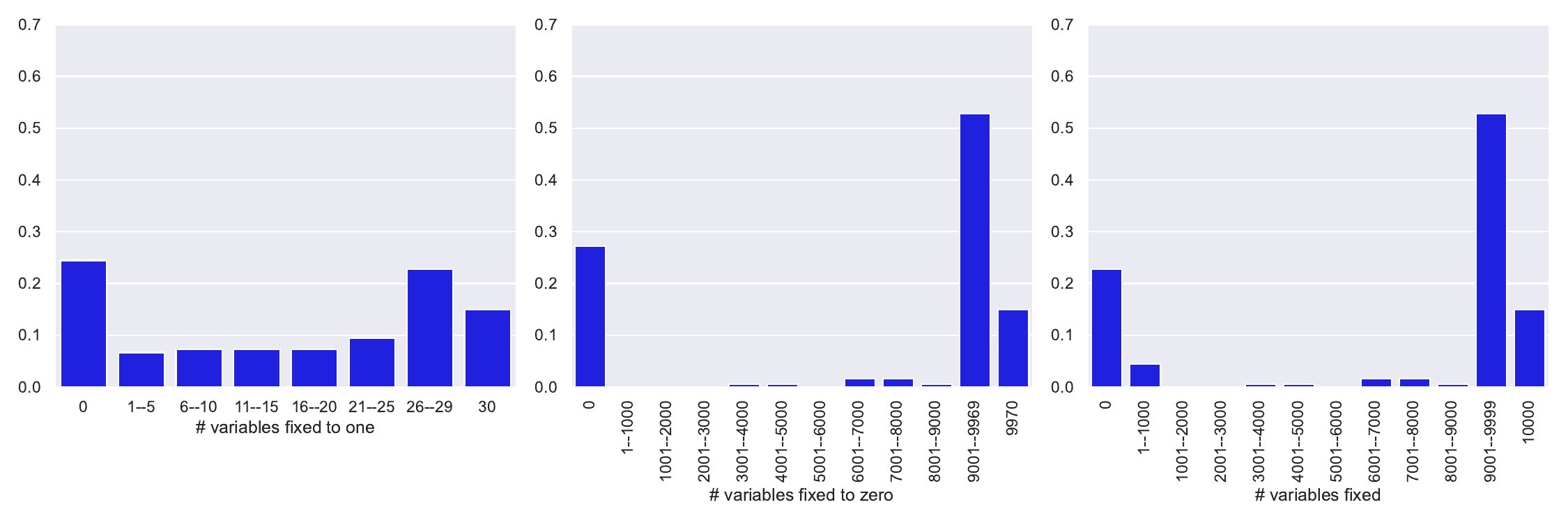}
  \caption{Distribution of the number of fixed variables.}
  \label{fig:ssr}
\end{figure}

Figure~\ref{fig:ssr} shows the distribution of the number of fixed variables (features) resulting from the safe screening across all synthetic datasets (180 instances).
The figure indicates that the safe screening successfully identifies the best subset of features in 15\% of the instances,
and more than 90\% of the features are fixed to zero or one in 67\% of the instances.
Conversely, in fewer than 23\% of the instances, no variables are fixed.

Next, we examine the effect of the safe screening in detail for each regime.
Table~\ref{tab:pegging} presents the results of our safe screening applied to the synthetic datasets.
In this table, the column labeled \# Fixed 1 (resp., \# Fixed 0) indicates the number of variables fixed to one (resp., zero),
and the column labeled Time represents the computation time required to execute Algorithm~\ref{alg:ssr} in seconds.
Each entry represents the average value across five trials, with the standard deviation shown in parentheses.
\begin{table}[h]
  \centering
  \caption{Results of safe screening.}
  \label{tab:pegging}
  \scalebox{0.9}{
  \begin{tabular}{ccccrrrrrr}
    \toprule
    $n$ & $\gamma$ & $\rho$ & $\sigma^2$ & \multicolumn{2}{c}{\textbf{\#~Fixed~1}} & \multicolumn{2}{c}{\textbf{\#~Fixed~0}} & \multicolumn{2}{c}{\textbf{Time (sec.)}} \\
    \toprule
    1000 & $2^{0}\gamma_{0}$ & 0.35 & 0.01 & 27.80 & (2.04) & 9969.40 & (0.54) & 10.89 & (0.53) \\
           &        &        & 0.10 & 28.20 & (0.83) & 9967.00 & (1.22) & 11.21 & (0.74) \\
           &        &        & 1.00 & 21.20 & (4.91) & 9955.80 & (7.22) & 10.30 & (0.84) \\
    \cmidrule{3-10}
           &        & 0.70 & 0.01 & 24.00 & (2.91) & 9964.20 & (1.78) & 10.008 & (0.39) \\
           &        &        & 0.10 & 25.00 & (2.54) & 9965.60 & (3.50) & 10.16 & (0.34) \\
           &        &        & 1.00 & 18.60 & (1.51) & 9951.20 & (12.02) & 9.36 & (0.44) \\
    \cmidrule{2-10}
           & $2^{2}\gamma_{0}$ & 0.35 & 0.01 & 25.60 & (2.88) & 9964.00 & (4.06) & 8.63 & (0.39) \\
           &        &        & 0.10 & 19.60 & (3.50) & 9932.00 & (21.11) & 8.49 & (0.51) \\
           &        &        & 1.00 & 3.00 & (2.34) & 2236.00 & (3328.96) & 7.92 & (0.08) \\
    \cmidrule{3-10}
           &        & 0.70 & 0.01 & 11.60 & (3.43) & 9918.40 & (11.41) & 7.74 & (0.22) \\
           &        &        & 0.10 & 8.20 & (4.60) & 9886.00 & (46.50) & 7.94 & (0.20) \\
           &        &        & 1.00 & 2.40 & (2.79) & 3256.60 & (4592.42) & 7.85 & (0.35) \\
    \cmidrule{2-10}
    & $2^{4}\gamma_{0}$ & 0.35 & 0.01 & 6.60 & (6.50) & 4690.20 & (4468.02) & 8.38 & (0.32) \\
    &        &        & 0.10 & 1.40 & (1.94) & 0.00 & (0.00) & 8.16 & (0.18) \\
    &        &        & 1.00 & 0.00 & (0.00) & 0.00 & (0.00) & 7.35 & ( 0.19) \\
\cmidrule{3-10}
    &        & 0.70 & 0.01 & 0.00 & (0.00) & 0.00 & (0.00) & 7.60 & (0.13) \\
    &        &        & 0.10 & 0.00 & (0.00) & 0.00 & (0.00) & 7.49 & (0.20) \\
    &        &        & 1.00 & 0.00 & (0.00) & 0.00 & (0.00) & 7.44 & (0.27) \\
    \midrule
    2000 & $2^{0}\gamma_{0}$ & 0.35 & 0.01 & 30.00 & (0.00) & 9970.00 & (0.00) & 28.54 & (0.85) \\
           &        &        & 0.10 & 29.80 & (0.44) & 9969.80 & (0.44) & 30.32 & (0.92) \\
           &        &        & 1.00 & 29.00 & (0.70) & 9968.40 & (1.14) & 28.40 & (1.55) \\
    \cmidrule{3-10}
           &        & 0.70 & 0.01 & 27.60 & (1.51) & 9966.60 & (1.34) & 26.40 & (1.20) \\
           &        &        & 0.10 & 27.40 & (1.34) & 9967.40 & (1.81) & 25.12 & (0.46) \\
           &        &        & 1.00 & 25.80 & (2.04) & 9965.80 & (2.28) & 25.80 & (0.85) \\
    \cmidrule{2-10}
           & $2^{2}\gamma_{0}$ & 0.35 & 0.01 & 29.80 & ( 0.44) & 9969.60 & (0.89) & 22.47 & (1.57) \\
           &        &        & 0.10 & 28.80 & (1.30) & 9969.00 & (1.00) & 23.16 & (0.92) \\
           &        &        & 1.00 & 24.60 & (2.79) & 9954.40 & (19.66) & 21.56 & (1.04) \\
    \cmidrule{3-10}
           &        & 0.70 & 0.01 & 19.60 & (5.94) & 9945.00 & (12.02) & 20.16 & (1.03) \\
           &        &        & 0.10 & 14.40 & (1.51) & 9946.20 & (7.79) & 19.88 & (0.87) \\
           &        &        & 1.00 & 9.20 & (1.48) & 9909.60 & (17.85) & 19.45 & (0.57) \\
    \cmidrule{2-10}
    & $2^{4}\gamma_{0}$ & 0.35 & 0.01 & 30.00 & (0.00) & 9970.00 & (0.00) & 18.33 & (0.78) \\
    &        &        & 0.10 & 28.80 & (2.68) & 9964.60 & (12.07) & 19.18 & (0.97) \\
    &        &        & 1.00 & 0.40 & (0.89) & 0.00 & (0.00) & 17.43 & (0.29) \\
\cmidrule{3-10}
    &        & 0.70 & 0.01 & 0.00 & (0.00) & 3076.00 & (4287.63) & 16.86 & (0.77) \\
    &        &        & 0.10 & 0.40 & (0.89) & 2327.00 & (3380.67) & 16.49 & (0.44) \\
    &        &        & 1.00 & 0.00 & (0.00) & 0.00 & (0.00) & 16.95 & (0.58) \\
    \bottomrule
    \end{tabular}
  }
\end{table}
From Table~\ref{tab:pegging}, we observe that the regularization parameter $\gamma$ has a significant impact on the performance.
As the parameter $\gamma$ decreases, which strengths the regularization term, the number of fixed variables increases.
We also observe that the correlation strength $\rho$ does not have a significant impact on the performance, while the noise level $\sigma^{2}$ plays a relatively important role.
In particular, when the noise level is low, the safe screening is more effective.
This trend is also observed in the experimental results of Atamturk and Gomez for ordinary regression problems.

Finally, to evaluate effect of the safe screening on the computational time reduction, we compare the following methods:
\begin{description}
  \item[BigM] solves the problem (Q) with the constraint~\eqref{eq:bigM} instead of the nonlinear constraint $w_{j}=w_{j}z_{j}$.
  \item[Persp] solves the perspective reformulation ($\text{P}_{\text{conic}}$).
  \item[Ours] solves the following problem reduced by Algorithm~\ref{alg:ssr}:
\end{description}
\[\begin{array}{lll}
  &\mini_{} & \displaystyle \frac{1}{n}\sum_{i=1}^{n} (t_{i}-y_{i}(\vt{w}^{\top}\vt{x}_{i}+b)+\log{y_{i}!}) +\frac{1}{\gamma} \sum_{j=1}^{m}s_{j}\\
  &\st & \displaystyle \sum_{j=1}^{m} z_{j} \le k,\\
  & & \displaystyle(t_{i}, 1, \vt{w}^{\top}\vt{x}_{i}+b) \in K_{\exp} \quad \text{for $i \in [n]$},\\
  & & \displaystyle (s_{j}/2, z_{j}, w_{j}) \in K_{\text{rq}}^{3} \quad \text{for $j \in [m]$},\\
  && \displaystyle z_{j}=0 \quad \text{for $j \in I_{0}$},\\
  && \displaystyle z_{j}=1 \quad \text{for $j \in I_{1}$},\\
  & & \displaystyle (\vt{w},b) \in \bb{R}^{m+1}, \vt{z} \in \{0,1\}^{m}, (\vt{s},\vt{t}) \in \bb{R}_{+}^{m+n}.
\end{array}
\]

In the BigM method, we set the parameter $M$ in the constraint~\eqref{eq:bigM} as $M=2$ since the true coefficient vector $\vt{w}^{\text{true}}$ is bounded by $1$.
Note that, in Ours, an upper bound solution obtained by executing the safe screening (Line~5 in Algorithm~\ref{alg:ssr}) is used as the initial solution in the Branch-and-Bound process.

The parameters in Mosek are set to the default values.
The computation of each methods is terminated if it does not finish within 600 seconds. In these cases, the results obtained within 600 seconds are taken as the final outcome.

Table~\ref{tab:bigM} (resp., Table~\ref{tab:persp} and Table~\ref{tab:ssr}) shows the computational performance of BigM (resp., Persp and Ours) for the synthetic datasets.
The column labeled Time represent the computation time in seconds,
and the column labeled Gap represent the optimality gap in percentage defined as
\begin{align*}
  \text{Gap} = \frac{\text{Obj}-\text{LB}}{\text{Obj}} \times 100,
\end{align*}
where $\text{Obj}$ is the objective value of the best feasible solution obtained by the method and $\text{LB}$ is the lower bound on the optimal value of the problem (P).
When the method cannot find any feasible solution within time limit, we set $\text{Gap}=100$.
The column labeled Nodes represent the number of nodes explored in the Branch-and-Bound tree, and the column labeled \text{OT} indicates the number of trials for which the method does not finish within 600 seconds.

Table~\ref{tab:bigM} demonstrates that BigM exhibits a substantial optimality gap (greater than $10\%$ in most cases) and requires exploring a significantly larger number of nodes in the Branch-and-Bound tree.

From Table~\ref{tab:persp}, it is shown that while the Persp method can optimally solve a few instances within the time limit, it fails to find even a feasible solution for most instances. This indicates the potential of perspective reformulation to provide tight lower bounds, as well as the inherent difficulty of mixed-integer conic optimization problems.

In contrast, Ours (Table~\ref{tab:ssr}) solves the vast majority of instances to optimality within 600 seconds. For the few remaining instances, the average optimality gap is less than approximately $1\%$,
and the number of nodes explored is significantly reduced compared to BigM.
Notably, even when the safe screening does not fix any variables (e.g., $n = 1000, \gamma = 2^{4}\gamma_{0}, \rho = 0.70$), our method still achieves a high-quality solution with an optimality gap of less than $0.80\%$ within the time limit. This performance is attributed to our greedy heuristic, which leverages the optimal dual solution $\bar{\vt{\lambda}}$ to provide a strong initial solution for the Branch-and-Bound process, leading to faster convergence.

\begin{table}
  \centering
  \caption{Results of BigM method.}
  \label{tab:bigM}
  \scalebox{0.9}{
  \begin{tabular}{cccc rrrrrrr}
    \toprule
    \multirow{2}{*}{$n$} & \multirow{2}{*}{$\gamma$} & \multirow{2}{*}{$\rho$} & \multirow{2}{*}{$\sigma^2$} &
    \multicolumn{7}{c}{\textbf{BigM}} \\
    \cmidrule{5-11}
    & & & & \multicolumn{2}{c}{\textbf{Time (sec.)}} & \multicolumn{2}{c}{\textbf{Gap (\%)}} & \multicolumn{2}{c}{\textbf{Nodes}} & \textbf{OT} \\
    \midrule
    1000 & $2^{0}\gamma_{0}$ & 0.35 & 0.01 & 600.00 & (0.00) & 11.29 & (0.64) & 698.80 & (95.24) & 5 \\
    & & & 0.10 & 600.00 & (0.00) & 11.88 & (0.58) & 662.20 & (95.71) & 5 \\
    & & & 1.00 & 600.00 & (0.00) & 16.15 & (0.79) & 732.60 & (79.74) & 5 \\
    \cmidrule{3-11}
    & & 0.70 & 0.01 & 600.00 & (0.00) & 11.00 & (0.36) & 662.00 & (95.88) & 5 \\
    & & & 0.10 & 600.00 & (0.00) & 11.83 & (0.54) & 664.00 & (94.53) & 5 \\
    & & & 1.00 & 600.00 & (0.00) & 16.94 & (0.67) & 664.20 & (93.84) & 5 \\
    \cmidrule{2-11}
    & $2^{2}\gamma_{0}$ & 0.35 & 0.01 & 600.00 & (0.00) & 22.19 & (0.87) & 732.80 & (76.47) & 5 \\
    & & & 0.10 & 600.00 & (0.00) & 23.19 & (1.17) & 697.20 & (96.57) & 5 \\
    & & & 1.00 & 600.00 & (0.00) & 31.88 & (0.90) & 664.20 & (93.84) & 5 \\
    \cmidrule{3-11}
    & & 0.70 & 0.01 & 600.00 & (0.00) & 21.35 & (0.93) & 626.00 & (79.45) & 5 \\
    & & & 0.10 & 600.00 & (0.00) & 23.84 & (0.75) & 733.20 & (76.13) & 5 \\
    & & & 1.00 & 600.00 & (0.00) & 31.75 & (0.79) & 664.80 & (93.30) & 5 \\
    \cmidrule{2-11}
    & $2^{4}\gamma_{0}$ & 0.35 & 0.01 & 600.00 & (0.00) & 27.64 & (0.65) & 733.40 & (76.81) & 5 \\
    & & & 0.10 & 600.00 & (0.00) & 30.01 & (0.89) & 768.60 & (2.07) & 5 \\
    & & & 1.00 & 600.00 & (0.00) & 42.50 & (0.62) & 733.00 & (76.02) & 5 \\
    \cmidrule{3-11}
    & & 0.70 & 0.01 & 600.00 & (0.00) & 28.10 & (0.61) & 700.60 & (95.97) & 5 \\
    & & & 0.10 & 600.00 & (0.00) & 29.91 & (0.42) & 769.40 & (2.50) & 5 \\
    & & & 1.00 & 600.00 & (0.00) & 42.25 & (1.00) & 769.00 & (2.73) & 5 \\
    \midrule
    2000 & $2^{0}\gamma_{0}$ & 0.35 & 0.01 & 600.00 & (0.00) & 4.69 & (0.22) & 99.00 & (2.23) & 5 \\
    & & & 0.10 & 600.00 & (0.00) & 4.89 & (0.15) & 100.00 & (0.00) & 5 \\
    & & & 1.00 & 600.00 & (0.00) & 7.21 & (0.18) & 92.20 & (5.21) & 5 \\
    \cmidrule{3-11}
    & & 0.70 & 0.01 & 600.00 & (0.00) & 5.33 & (0.36) & 97.80 & (2.58) & 5 \\
    & & & 0.10 & 600.00 & (0.00) & 5.85 & (0.18) & 100.00 & (0.00) & 5 \\
    & & & 1.00 & 600.00 & (0.00) & 7.60 & (0.19) & 89.00 & (5.04) & 5 \\
    \cmidrule{2-11}
    & $2^{2}\gamma_{0}$ & 0.35 & 0.01 & 600.00 & (0.00) & 12.83 & (0.38) & 89.00 & (8.09) & 5 \\
    & & & 0.10 & 600.00 & (0.00) & 13.97 & (0.80) & 93.40 & (2.19) & 5 \\
    & & & 1.00 & 600.00 & (0.00) & 20.01 & (0.57) & 95.80 & (1.09) & 5 \\
    \cmidrule{3-11}
    & & 0.70 & 0.01 & 600.00 & (0.00) & 15.02 & (0.46) & 93.60 & (1.94) & 5 \\
    & & & 0.10 & 600.00 & (0.00) & 14.59 & (0.59) & 92.80 & (3.19) & 5 \\
    & & & 1.00 & 600.00 & (0.00) & 20.00 & (0.61) & 94.60 & (0.89) & 5 \\
    \cmidrule{2-11}
    & $2^{4}\gamma_{0}$ & 0.35 & 0.01 & 600.00 & (0.00) & 21.18 & (0.27) & 97.60 & (3.43) & 5 \\
    & & & 0.10 & 600.00 & (0.00) & 23.39 & (0.51) & 95.00 & (4.52) & 5 \\
    & & & 1.00 & 600.00 & (0.00) & 33.51 & (0.78) & 95.00 & (5.19) & 5 \\
    \cmidrule{3-11}
    & & 0.70 & 0.01 & 600.00 & (0.00) & 22.15 & (0.75) & 93.60 & (6.98) & 5 \\
    & & & 0.10 & 600.00 & (0.00) & 24.20 & (0.76) & 96.20 & (3.56) & 5 \\
    & & & 1.00 & 600.00 & (0.00) & 33.26 & (0.23) & 98.60 & (3.13) & 5 \\
    \bottomrule
  \end{tabular}
  }
\end{table}

\begin{table}
  \centering
  \caption{Results of Persp method.}
  \label{tab:persp}
  \scalebox{0.9}{
  \begin{tabular}{cccc rrrrrrr}
    \toprule
    \multirow{2}{*}{$n$} & \multirow{2}{*}{$\gamma$} & \multirow{2}{*}{$\rho$} & \multirow{2}{*}{$\sigma^2$} &
    \multicolumn{7}{c}{\textbf{Persp}} \\
    \cmidrule{5-11}
    & & & & \multicolumn{2}{c}{\textbf{Time (sec.)}} & \multicolumn{2}{c}{\textbf{Gap (\%)}} & \multicolumn{2}{c}{\textbf{Nodes}} & \textbf{OT} \\
    \midrule
    1000 & $2^{0}\gamma_{0}$ & 0.35 & 0.01 & 504.67 & (213.73) & 80.00 & (44.72) & 0.00 & (0.00) &4 \\
    & & & 0.10 & 600.00 & (0.00) & 100.00 & (0.00) & 0.00 & (0.00) & 5 \\
    & & & 1.00 & 600.00 & (0.00) & 100.00 & (0.00) & 0.00 & (0.00) & 5\\
    \cmidrule{3-11}
    & & 0.70 & 0.01 & 600.00 & (0.00) & 100.00 & (0.00) & 0.00 & (0.00) & 5 \\
    & & & 0.10 & 600.00 & (0.00) & 100.00 & (0.00) & 0.00 & (0.00) & 5 \\
    & & & 1.00 & 600.00 & (0.00) & 100.00 & (0.00) & 0.00 & (0.00) & 5 \\
    \cmidrule{2-11}
    & $2^{2}\gamma_{0}$ & 0.35 & 0.01 & 600.00 & (0.00) & 100.00 & (0.00) & 0.00 & (0.00) & 5 \\
    & & & 0.10 & 600.00 & (0.00) & 100.00 & (0.00) & 0.00 & (0.00) & 5 \\
    & & & 1.00 & 600.00 & (0.00) & 100.00 & (0.00) & 0.00 & (0.00) & 5 \\
    \cmidrule{3-11}
    & & 0.70 & 0.01 & 600.00 & (0.00) & 100.00 & (0.00) & 0.00 & (0.00) & 5 \\
    & & & 0.10 & 600.00 & (0.00) & 100.00 & (0.00) & 0.00 & (0.00) & 5\\
    & & & 1.00 & 600.00 & (0.00) & 100.00 & (0.00) & 0.00 & (0.00) & 5\\
    \cmidrule{2-11}
    & $2^{4}\gamma_{0}$ & 0.35 & 0.01 & 600.00 & (0.00) & 100.00 & (0.00) & 0.00 & (0.00) &5 \\
    & & & 0.10 & 600.00 & (0.00) & 100.00 & (0.00) & 0.00 & (0.00) &5 \\
    & & & 1.00 & 600.00 & (0.00) & 1.19 & (0.53) & 46.40 & (13.63) &5 \\
    \cmidrule{3-11}
    & & 0.70 & 0.01 & 600.00 & (0.00) & 100.00 & (0.00) & 0.00 & (0.00) &5 \\
    & & & 0.10 & 600.00 & (0.00) & 100.00 & (0.00) & 0.00 & (0.00) & 5\\
    & & & 1.00 & 600.00 & (0.00) & 1.01 & (0.48) & 41.00 & (18.28) &5 \\
    \midrule
    2000 & $2^{0}\gamma_{0}$ & 0.35 & 0.01 & 600.00 & (0.00) & 100.00 & (0.00) & 0.00 & (0.00) &5 \\
    & & & 0.10 & 496.25 & (232.96) & 80.00 & (44.72) & 0.00 & (0.00) &4 \\
    & & & 1.00 & 600.00 & (0.00) & 100.00 & (0.00) & 0.00 & (0.00) &5\\
    \cmidrule{3-11}
    & & 0.70 & 0.01 & 519.76 & (180.33) & 80.00 & (44.72) & 0.00 & (0.00) &4 \\
    & & & 0.10 & 600.00 & (0.00) & 100.00 & (0.00) & 0.00 & (0.00) &5\\
    & & & 1.00 & 519.06 & (181.70) & 80.00 & (44.72) & 0.00 & (0.00) & 4\\
    \cmidrule{2-11}
    & $2^{2}\gamma_{0}$ & 0.35 & 0.01 & 296.85 & (278.47) & 40.00 & (54.77) & 0.00 & (0.00) &2 \\
    & & & 0.10 & 388.83 & (289.65) &60.00 & (54.77) & 0.00 & (0.00) & 3\\
    & & & 1.00 & 600.00 & (0.00) & 100.00 & (0.00) & 0.00 & (0.00) &5\\
    \cmidrule{3-11}
    & & 0.70 & 0.01 & 600.00 & (0.00) & 100.00 & (0.00) & 0.00 & (0.00) &5\\
    & & & 0.10 & 600.00 & (0.00) & 100.00 & (0.00) & 0.00 & (0.00) & 5\\
    & & & 1.00 & 600.00 & (0.00) & 100.00 & (0.00) & 0.00 & (0.00) & 5\\
    \cmidrule{2-11}
    & $2^{4}\gamma_{0}$ & 0.35 & 0.01 & 174.78 & (238.03) & 20.00 & (44.72) & 0.00 & (0.00) & 1 \\
    & & & 0.10 & 176.75 & (236.72) & 20.00 & (44.72) & 0.00 & (0.00) &1\\
    & & & 1.00 & 600.00 & (0.00) & 100.00 & (0.00) & 0.00 & (0.00) &5\\
    \cmidrule{3-11}
    & & 0.70 & 0.01 & 600.00 & (0.00) & 100.00 & (0.00) & 0.00 & (0.00) &5 \\
    & & & 0.10 & 600.00 & (0.00) & 100.00 & (0.00) & 0.00 & (0.00) & 5\\
    & & & 1.00 & 600.00 & (0.00) & 100.00 & (0.00) & 0.00 & (0.00) &5\\
    \bottomrule
  \end{tabular}
  }
\end{table}

\begin{table}
  \centering
  \caption{Results of our method (Ours).}
  \label{tab:ssr}
  \scalebox{0.9}{
  \begin{tabular}{cccc rrrrrr}
    \toprule
    \multirow{2}{*}{$n$} & \multirow{2}{*}{$\gamma$} & \multirow{2}{*}{$\rho$} & \multirow{2}{*}{$\sigma^2$} &
    \multicolumn{6}{c}{\textbf{Ours}}\\
    \cmidrule{5-10}
    & & & & \multicolumn{2}{c}{\textbf{Time (sec.)}} & \multicolumn{2}{c}{\textbf{Gap (\%)}} & \multicolumn{2}{c}{\textbf{Nodes}} \\
    \midrule
    1000 & $2^{0}\gamma_{0}$ & 0.35 & 0.01 & 30.64 & (5.11) & 0.00 & (0.00) & 0.00 & (0.00) \\
    & & & 0.10 & 26.88 & (2.24) & 0.00 & (0.00) & 0.00 & (0.00) \\
    & & & 1.00 & 34.47 & (9.83) & 0.00 & (0.00) & 0.00 & (0.00) \\
    \cmidrule{3-10}
    & & 0.70 & 0.01 & 26.25 & (0.60) & 0.00 & (0.00) & 0.00 & (0.00) \\
    & & & 0.10 & 26.92 & (2.27) & 0.00 & (0.00) & 0.00 & (0.00) \\
    & & & 1.00 & 34.79 & (10.68) & 0.00 & (0.00) & 0.00 & (0.00) \\
    \cmidrule{2-10}
    & $2^{2}\gamma_{0}$ & 0.35 & 0.01 & 41.89 & (19.06) & 0.00 & (0.00) & 0.60 & (1.34) \\
    & & & 0.10 & 90.23 & (31.90) & 0.00 & (0.00) & 10.00 & (7.71) \\
    & & & 1.00 & 356.34 & (224.56) & 0.01 & (0.02) & 113.60 & (100.22) \\
    \cmidrule{3-10}
    & & 0.70 & 0.01 & 108.96 & (65.20) & 0.00 & (0.00) & 48.20 & (81.49) \\
    & & & 0.10 & 143.95 & (68.25) & 0.00 & (0.00) & 59.40 & (66.94) \\
    & & & 1.00 & 499.14 & (165.41) & 0.01 & (0.01) & 195.60 & (88.09) \\
    \cmidrule{2-10}
    & $2^{4}\gamma_{0}$ & 0.35 & 0.01 & 194.55 & (198.48) & 0.00 & (0.00) & 62.60 & (121.19) \\
    & & & 0.10 & 324.69 & (83.43) & 0.00 & (0.00) & 94.20 & (50.8) \\
    & & & 1.00 & 601.08 & (0.24) & 1.02 & (0.44) & 250.20 & (77.82) \\
    \cmidrule{3-10}
    & & 0.70 & 0.01 & 391.14 & (167.39) & 0.04 & (0.09) & 122.40 & (91.77) \\
    & & & 0.10 & 505.28 & (133.87) & 0.06 & (0.06) & 200.80 & (100.11) \\
    & & & 1.00 & 600.99 & (0.24) & 0.80 & (0.29) & 248.40 & (77.98) \\
    \midrule
    2000 & $2^{0}\gamma_{0}$ & 0.35 & 0.01 & 61.55 & (0.13) & 0.00 & (0.00) & 0.00 & (0.00) \\
    & & & 0.10 & 61.88 & (0.53) & 0.00 & (0.00) & 0.00 & (0.00) \\
    & & & 1.00 & 63.14 & (1.01) & 0.00 & (0.00) & 0.00 & (0.00) \\
    \cmidrule{3-10}
    & & 0.70 & 0.01 & 62.21 & (0.46) & 0.00 & (0.00) & 0.00 & (0.00) \\
    & & & 0.10 & 62.60 & (0.83) & 0.00 & (0.00) & 0.00 & (0.00) \\
    & & & 1.00 & 67.51 & (9.17) & 0.00 & (0.00) & 0.00 & (0.00) \\
    \cmidrule{2-10}
    & $2^{2}\gamma_{0}$ & 0.35 & 0.01 & 62.92 & (0.33) & 0.00 & (0.00) & 0.00 & (0.00) \\
    & & & 0.10 & 62.95 & (0.35) & 0.00 & (0.00) & 0.00 & (0.00) \\
    & & & 1.00 & 71.03 & (18.11) & 0.00 & (0.00) & 0.00 & (0.00) \\
    \cmidrule{3-10}
    & & 0.70 & 0.01 & 135.78 & (79.35) & 0.00 & (0.00) & 2.60 & (5.81) \\
    & & & 0.10 & 181.62 & (69.99) & 0.00 & (0.00) & 6.00 & (5.48) \\
    & & & 1.00 & 209.56 & (96.91) & 0.00 & (0.00) & 8.80 & (8.29) \\
    \cmidrule{2-10}
    & $2^{4}\gamma_{0}$ & 0.35 & 0.01 & 61.02 & (0.3) & 0.00 & (0.00) & 0.00 & (0.00) \\
    & & & 0.10 & 70.48 & (20.17) & 0.00 & (0.00) & 0.00 & (0.00) \\
    & & & 1.00 & 494.45 & (146.72) & 0.04 & (0.06) & 36.40 & (19.67) \\
    \cmidrule{3-10}
    & & 0.70 & 0.01 & 535.86 & (75.20) & 0.01 & (0.01) & 60.80 & (34.40) \\
    & & & 0.10 & 569.48 & (71.18) & 0.02 & (0.02) & 61.40 & (13.16) \\
    & & & 1.00 & 602.23 & (0.72) & 0.16 & (0.14) & 52.40 & (8.17) \\
    \bottomrule
  \end{tabular}
  }
\end{table}

\subsection{Results for higher dimensional datasets}

To further assess the scalability of our method, we conduct experiments utilizing higer dimensional datasets with $m\in \{30000, 50000\}$ and $n\in \{1000, 2000\}$.
Data generation followed the methodology outlined in Section 4.1; however, to focus on challenging conditions, experiments are limited to scenarios exhibiting high feature correlation and noise levels $(\rho,\sigma^{2})=(0.70, 1.00)$.
The parameter $k$ is set to $30$, and the time limit for each instance is set to 1800 seconds.

Table~\ref{tab:highdim} presents the results of our method, which demonstrates its ability to effectively handle larger datasets.
In Table~\ref{tab:highdim}, the columns labeled \#~Fixed indicate the number of fixed variables, while Time~1 (resp., Time~2) represents the computation time required to execute Algorithm~\ref{alg:ssr} (resp., required to solve the reduced problem). The other columns are defined similarly to those in Table~\ref{tab:ssr}.

Table~\ref{tab:highdim} shows the scalability of our method with higher-dimensional datasets ($m \in \{30000, 50000\}$ and $n \in \{1000, 2000\}$).
While increasing the regularization parameter $\gamma$ sometimes leads to reaching the time limit (1800 seconds), the method consistently finds solutions with a remarkably small optimality gap.
Specifically, with the lowest regularization parameter $2^{0}\gamma_{0}$,
the method often fixes a large number of variables and quickly solves the problem  to optimality.
Even when higher regularization parameters $2^{2}\gamma_{0}$ and $2^{4}\gamma_{0}$ result in reaching the time limit, the solutions obtained still exhibit very small optimality gaps, demonstrating the robustness and effectiveness of our approach in finding near-optimal solutions for high-dimensional datasets.

However, for the largest instances ($m = 50000$ and $n = 2000$), memory overflow errors occurred in one out of five trials each with the regularization parameters $2^{2}\gamma_{0}$ and $2^{4}\gamma_{0}$, indicating a potential limitation of the current implementation for extremely large datasets.

\begin{table}[h]
  \centering
  \caption{Results for higher dimensional datasets.}
  \label{tab:highdim}
  \scalebox{0.8}{
  \begin{tabular}{cccrrrrrrrrrrr}
    \toprule
    $m$ & $n$ & $\gamma$ & \multicolumn{2}{c}{\textbf{\#~Fixed}} & \multicolumn{2}{c}{\textbf{Time~1 (sec.)}} & \multicolumn{2}{c}{\textbf{Time~2 (sec.)}} & \multicolumn{2}{c}{\textbf{Gap (\%)}} & \multicolumn{2}{c}{\textbf{Nodes}} & \textbf{OT} \\
    \midrule
    30000 & 1000 & $2^{0}\gamma_{0}$ & 29965.00 & (11.79) & 48.23 & (22.46) & 120.25 & (18.26) & 0.00 & (0.00) & 0.00 & (0.00) & 0 \\
    & & $2^{2}\gamma_{0}$ & 0.20 & (0.45) & 51.67 & (43.61) & 1800.00 & (0.00) & 0.05 & (0.04) & 103.00 & (29.04) & 5 \\
    & & $2^{4}\gamma_{0}$ & 0.00 & (0.00) & 54.32 & (57.36) & 1800.00 & (0.00) & 1.24 & (0.24) & 83.00 & (23.03) & 5 \\
    \cmidrule{2-14}
    & 2000 & $2^{0}\gamma_{0}$ & 29993.80 & (1.30) & 93.20 & (8.59) & 307.12 & (40.17) & 0.00 & (0.00) & 0.00 & (0.00) & 0 \\
    & & $2^{2}\gamma_{0}$ & 29387.00 & (1045.20) & 72.03 & (2.87) & 759.30 & (535.91) & 0.00 & (0.00) &  9.00 & (14.04) & 0 \\
    & & $2^{4}\gamma_{0}$ &  0.00 & (0.00) & 59.25 & (1.86) & 1800.00 & (0.00) & 0.10 & (0.06) & 32.20 & (2.95) & 5 \\
    \midrule
    50000 & 1000 & $2^{0}\gamma_{0}$ & 49943.80 & (26.86) & 207.47 & (173.69) & 271.02 & (34.39) & 0.00 & (0.00) & 0.00 & (0.00) & 0 \\
    &  & $2^{2}\gamma_{0}$ & 0.40 & (0.55) & 227.01 & (156.52) & 1800.00 & (0.00) & 0.06 & (0.06) & 2.00 & (1.22) & 5\\
    &  & $2^{4}\gamma_{0}$ & 0.00 & (0.00) & 105.57 & (80.02) & 1800.00 & (0.00) & 1.25 & (0.83) & 1.80 & (1.30) & 5\\
    \cmidrule{2-14}
    & 2000 & $2^{0}\gamma_{0}$ & 49992.80 & (3.42) & 177.97 & (37.98) & 592.90 & (117.47) & 0.00 & (0.00) & 0.00 & (0.00) & 0 \\
    &  & $2^{2}\gamma_{0}$ & 49756.75 & (36.12) & 158.55 & (7.74) & 1310.85 & (582.83) & 0.01 & (0.01) & 5.00 & (8.72) & 2 \\
    &  & $2^{4}\gamma_{0}$ & 0.00 & (0.00) & 114.59 & (2.91) & 1800.00 & (0.00) & 0.22 & (0.38) & 0.00 & (0.00) & 4 \\
    \bottomrule
  \end{tabular}
  }
\end{table}

\section{Conclusion}

In this paper, we propose an approach to cardinality-constrained Poisson regression for high-dimensional count data. We formulate the problem as a mixed-integer conic optimization problem and employ a safe screening based on Fenchel conjugates to enhance computational efficiency.

Our computational experiments on synthetic datasets demonstrate that the safe screening significantly reduces problem size, leading to substantial improvements in computational time. Specifically, our approach solves instances with up to 50,000 candidate features and 2,000 observations within a reasonable time limit (1,800 seconds). For example, with low regularization parameters and a high noise situation $(\rho,\sigma^{2})=(0.70, 1.00)$, we achieved optimality for instances with 50,000 features by fixing nearly all variables using the safe screening.
Although the safe screening is less effective with high regularization parameters, it provides a good starting point for the optimization process through greedy heuristics, resulting in high-quality solutions with optimality gaps below $1.25\%$ in most cases.

A future direction is further investigation of the screening rule to enhance its effectiveness, especially in scenarios with high noise levels.
While our experiments on synthetic datasets indicate the method's efficacy for high-dimensional data, additional experiments on real-world datasets are necessary to validate its practical applicability.
Another direction is to extend the proposed method to more general regression models, such as the negative binomial regression model, which is widely used across various fields.

\section*{Acknowledgement}
This work is supported by JSPS Grant-in-Aid (20KK0106 and 22K17856).

\bibliographystyle{abbrvnat}
\bibliography{ref}

\begin{thebibliography}{28}
\providecommand{\natexlab}[1]{#1}
\providecommand{\url}[1]{\texttt{#1}}
\expandafter\ifx\csname urlstyle\endcsname\relax
  \providecommand{\doi}[1]{doi: #1}\else
  \providecommand{\doi}{doi: \begingroup \urlstyle{rm}\Url}\fi

\bibitem[Ahari and Kocuk(2023)]{AK_2023}
S.~A. Ahari and B.~Kocuk.
\newblock A mixed-integer exponential cone programming formulation for feature
  subset selection in logistic regression.
\newblock \emph{EURO Journal on Computational Optimization}, 11:\penalty0
  100069, 2023.

\bibitem[Algamal(2019)]{A_SOIC_2019}
Z.~Y. Algamal.
\newblock Variable selection in count data regression model based on firefly
  algorithm.
\newblock \emph{Statistics, Optimization and Information Computing},
  7:\penalty0 520--529, 2019.

\bibitem[Andres and Huber(2010)]{AH_GB2010}
S.~Andres and W.~Huber.
\newblock Differential expression analysis for sequence count data.
\newblock \emph{Genome Biology}, 11:\penalty0 R106, 2010.
\newblock \doi{10.1186/gb-2010-11-10-r106}.

\bibitem[Atamt\"{u}rk and G\'{o}mez(2020)]{AG_ICML2020}
A.~Atamt\"{u}rk and A.~G\'{o}mez.
\newblock Safe screening rules for $\ell_{0}$-regression from perspective
  relaxations.
\newblock In \emph{Proceedings of the 37th International Conference on Machine
  Learning}, pages 421--430. PMLR, 2020.

\bibitem[Balas and Carrera(1996)]{BC_OR1996}
E.~Balas and M.~C. Carrera.
\newblock A dynamic subgradient-based branch-and-bound procedure for set
  covering.
\newblock \emph{Operations research}, 44:\penalty0 875--890, 1996.

\bibitem[Bertsimas et~al.(2016)Bertsimas, King, and Mazumder]{BKM_AOR2016}
D.~Bertsimas, A.~King, and R.~Mazumder.
\newblock Best subset selection via a modern optimization lens.
\newblock \emph{The Annals of Statistics}, 44:\penalty0 813--852, 2016.

\bibitem[Bertsimas et~al.(2020)Bertsimas, Pauphilet, and Parys]{BPP_2020}
D.~Bertsimas, J.~Pauphilet, and B.~V. Parys.
\newblock Sparse regression: Scalable algorithms and empirical performance.
\newblock \emph{Statistical Science}, 35:\penalty0 555--578, 2020.

\bibitem[Bonnefoy et~al.(2014)Bonnefoy, Emiya, Ralaivola, and
  Gribonval]{BERG_EUSIPCO2014}
A.~Bonnefoy, V.~Emiya, L.~Ralaivola, and R.~Gribonval.
\newblock A dynamic screening principle for the lasso.
\newblock In \emph{22nd European Signal Processing Conference}, pages 6--10,
  2014.

\bibitem[Dedieu et~al.(2021)Dedieu, Hazimeh, and Mazumder]{DHHM_JMLR2021}
A.~Dedieu, H.~Hazimeh, and R.~Mazumder.
\newblock Learning sparse classifiers: Continuous and mixed integer
  optimization perspectives.
\newblock \emph{Journal of Machine Learning Research}, 22:\penalty0 1--47,
  2021.

\bibitem[Efroymson(1960)]{E_MMDC1960}
M.~A. Efroymson.
\newblock \emph{Multiple Regression Analysis}.
\newblock John Wiley \& Sons, 1960.

\bibitem[Fitzgerald et~al.(2022)Fitzgerald, Jones, and Engelhardt]{FJE_BMC2022}
T.~Fitzgerald, A.~Jones, and B.~E. Engelhardt.
\newblock A poisson reduced-rank regression model for association mapping in
  sequencing data.
\newblock \emph{BMC Bioinformatics}, 23:\penalty0 529, 2022.

\bibitem[Forman(2003)]{F_JMLR2003}
G.~Forman.
\newblock An extensive empirical study of feature selection metrics for text
  classification.
\newblock \emph{Journal of Machine Learning Research}, 3:\penalty0 1289--1305,
  2003.

\bibitem[G{\'o}mez and Prokopyev(2021)]{GP_INFORMS2021}
A.~G{\'o}mez and O.~A. Prokopyev.
\newblock A mixed-integer fractional optimization approach to best subset
  selection.
\newblock \emph{INFORMS Journal on Computing}, 33:\penalty0 551--565, 2021.

\bibitem[G{\"u}nl{\"u}k and Linderoth(2010)]{GL_MP2010}
O.~G{\"u}nl{\"u}k and J.~Linderoth.
\newblock Perspective reformulations of mixed integer nonlinear programs with
  indicator variables.
\newblock \emph{Mathematical programming}, 124:\penalty0 183--205, 2010.

\bibitem[Hastie et~al.(2020)Hastie, Tibshirani, and Tibshirani]{HTT_STS2020}
T.~Hastie, R.~Tibshirani, and R.~Tibshirani.
\newblock Best subset, forward stepwise or lasso? {A}nalysis and
  recommendations based on extensive comparisons.
\newblock \emph{Statistical Science}, 35:\penalty0 579--592, 2020.

\bibitem[Hiriart-Urruty and Lemar{\'e}cha(1993)]{HU_2011}
J.-B. Hiriart-Urruty and C.~Lemar{\'e}cha.
\newblock \emph{Convex {A}nalysis and {M}inimization {A}lgorithms I:
  {F}undamentals}.
\newblock Springer, 1st edition, 1993.

\bibitem[Ko{\c{c}} et~al.(2018)Ko{\c{c}}, D{\"u}nder, G{\"u}m{\"u}{\c{s}}tekin,
  Ko{\c{c}}, and Cengiz]{AK_CSTM2018}
H.~Ko{\c{c}}, E.~D{\"u}nder, S.~G{\"u}m{\"u}{\c{s}}tekin, T.~Ko{\c{c}}, and
  M.~A. Cengiz.
\newblock Particle swarm optimization-based variable selection in poisson
  regression analysis via information complexity-type criteria.
\newblock \emph{Communications in Statistics Thoery and Methods}, 47:\penalty0
  5298--5306, 2018.

\bibitem[Kocuk(2021)]{K_IJOC2021}
B.~Kocuk.
\newblock Conic reformulations for {K}ullback-{L}eibler divergence constrained
  distributionally robust optimization and applications.
\newblock \emph{An International Journal of Optimization and Control: {T}heory
  \& {A}pplications}, 11:\penalty0 139--151, 2021.

\bibitem[Miller(2002)]{A_2002}
A.~Miller.
\newblock \emph{Subset {S}election in {R}egression}.
\newblock Taylor \& Francis, 2nd edition, 2002.

\bibitem[Miyashiro and Takano(2015)]{MT_EJOR2015}
R.~Miyashiro and Y.~Takano.
\newblock Mixed integer second-order cone programming formulations for variable
  selection in linear regression.
\newblock \emph{European Journal of Operational Research}, 247:\penalty0
  721--731, 2015.

\bibitem[{MOSEK ApS}(2024{\natexlab{a}})]{Mosek_11}
{MOSEK ApS}.
\newblock {MOSEK} {O}ptimizer {API} for {J}ulia 11.0.13.
\newblock 2024{\natexlab{a}}.
\newblock URL \url{https://docs.mosek.com/11.0/toolbox/julia/index.html}.

\bibitem[{MOSEK ApS}(2024{\natexlab{b}})]{Mosek_cookbook}
{MOSEK ApS}.
\newblock {MOSEK} {M}odeling {C}ookbook~3.3.0, 2024{\natexlab{b}}.
\newblock URL \url{https://docs.mosek.com/modeling-cookbook/}.

\bibitem[Natarajan(1995)]{N_SICOMP1995}
B.~K. Natarajan.
\newblock Sparse approximate solutions to linear systems.
\newblock \emph{SIAM Journal on Computing}, 24:\penalty0 227--234, 1995.

\bibitem[Nauss(1976)]{N_MS1976}
R.~M. Nauss.
\newblock An efficient algorithm for the 0-1 knapsack problem.
\newblock \emph{Management Science}, 23:\penalty0 27--31, 1976.

\bibitem[Saishu et~al.(2021)Saishu, Kudo, and Takano]{SKT_PLOSONE2021}
H.~Saishu, K.~Kudo, and Y.~Takano.
\newblock Sparse poisson regression via mixed-integer optimization.
\newblock \emph{PLOS ONE}, 16:\penalty0 e0249916, 2021.
\newblock \doi{https://doi.org/10.1371/journal.pone.0248536}.

\bibitem[Sato et~al.(2016)Sato, Takano, Miyashiro, and Yoshise]{STMY_COAP2016}
T.~Sato, Y.~Takano, R.~Miyashiro, and A.~Yoshise.
\newblock Feature subset selection for logistic regression via mixed integer
  optimization.
\newblock \emph{Computational Optimization and Applications}, 64:\penalty0
  865--880, 2016.

\bibitem[Townes et~al.(2019)Townes, Hicks, Aryee, and Irizarry]{THAI_GB2019}
F.~W. Townes, S.~C. Hicks, M.~J. Aryee, and R.~A. Irizarry.
\newblock Feature selection and dimension reduction for single-cell rna-seq
  based on a multinomial model.
\newblock \emph{Genome Biology}, 20:\penalty0 1--16, 2019.

\bibitem[Ye and Xie(2021)]{YX_arxiv2021}
Q.~Ye and W.~Xie.
\newblock Second-order conic and polyhedral approximations of the exponential
  cone: application to mixed-integer exponential conic programs.
\newblock \emph{arXiv}, 2106.09123, 2021.

\end{thebibliography}

\end{document}